\documentclass[11pt]{article}
\usepackage{amsthm}
\usepackage{amsmath}  
\usepackage{amssymb}  
\usepackage{geometry}  
\usepackage{hyperref}  
\usepackage{authblk}  
\usepackage{setspace}  
\usepackage{pxfonts}
\usepackage{yfonts}
\usepackage{dsfont}
\usepackage{graphicx}
\usepackage{relsize}
\usepackage{color}
\usepackage{marvosym}
\usepackage{graphicx}
\usepackage{relsize}
\usepackage[latin1]{inputenc}
\usepackage{enumerate}

\def\bel{\begin{equation}\label}
\def\eeq{\end{equation}}
\def\ds{\displaystyle}

\def\mt{\longrightarrow}

\def\Rec{{\bf R}}
\def\R{\mathbb R}

\def\C{\mathfrak{B}}

\def\N{{\bf N}}

\def\A{{\bf A}}

\def\L{{\bf L}}
\def\U{{\bf U}}

\def\Hat{\widehat}

\def\vol{{\bf vol}}

\def\M{{\bf M}}

\def\Cup{{\bigcup}}

\def\alpha{\alphaup}
\def\beta{\betaup}
\def\gamma{\gammaup}
\def\delta{\deltaup}
\def\theta{\thetaup}
\def\xi{{\xiup}}
\def\eta{{\etaup}}
\def\tau{{\tauup}}
\def\rho{{\rhoup}}
\def\phi{{\phiup}}
\def\psi{{\psiup}}
\def\lambda{{\lambdaup}}
\def\omega{\omegaup}
\def\varphi{{\varphiup}}
\def\gamma{{\gammaup}}

\newtheorem{remark}{Remark}[section]

\begin{document}

\[\begin{array}{cc}\ds\hbox{\LARGE{\bf On a family of strong fractional maximal operators}}
\end{array}\]

\[\hbox{Zipeng Wang}\]
\begin{abstract}
We study the strong fractional maximal function defined as
\[\M_\rho^\alpha f(x)~=~\sup_{\Rec\ni x}~ \vol\{\Rec\}^{\alpha-1}\int_\Rec |f(y_1+\rho_1,y_2+\rho_2,\ldots,y_{\N-1}+\rho_{\N-1},y_\N)|dy,\qquad 0\leq\alpha<1\]
where $\Rec\subset\R^\N$ is a rectangle parallel to the coordinates. Moreover, 
\[\rho_i~=~\rho_i(y_{i+1},\ldots,y_\N,x_{i+1},\ldots,x_\N),\qquad i=1,2,\ldots,\N-1\]
are measurable functions. We prove $\M_\rho^\alpha\colon\L^p(\R^\N)\mt\L^q(\R^\N)$  for $\alpha={1\over p}-{1\over q},1<p\leq q<\infty$.
\end{abstract}

\section{Introduction}
\setcounter{equation}{0}
Strong maximal functions play an important role in the multi-parameter theory of harmonic analysis.
Consider 
\bel{M_rho}
\M_\rho^\alpha f(x)~=~\sup_{\Rec\ni x}~ \vol\{\Rec\}^{\alpha-1}\int_\Rec |f(y_1+\rho_1,y_2+\rho_2,\ldots,y_{\N-1}+\rho_{\N-1},y_\N)|dy,\qquad 0\leq\alpha<1
\eeq
where $\Rec\subset\R^\N$ is a rectangle parallel to the coordinates. Moreover, 
\bel{rho_i}\rho_i~=~\rho_i(y_{i+1},\ldots,y_\N,x_{i+1},\ldots,x_\N),\qquad i=1,2,\ldots,\N-1
\eeq
are some measurable functions. A particular example would be the strong maximal operator defined on Heisenberg groups. Let 
\[
\M_\mu f(u,v,t)~=~\sup_{\Rec\ni(0,0,0)}
\vol \{\Rec\}^{-1} \iiint_\Rec\left|f [(u,v,t)\odot(\xi,\eta,\tau)^{-1}]\right|d\xi d\eta d\tau
\]
where $\odot$ is  a  multiplication law:
$(u,v,t)\odot(\xi,\eta,\tau)=\big[u+\xi, v+\eta,t+\tau+\mu(u\cdot\eta-v\cdot\xi)\big],~\mu\neq0$
for every $(u,v,t)\in\R^n\times\R^n\times\R$ and $(\xi,\eta,\tau)^{-1}=(-\xi,-\eta,-\tau)\in\R^n\times\R^n\times\R$.

By taking $\xi\mt u-\xi$, $ \eta\mt v-\eta$ and $\tau\mt t-\tau$, we find $\M_\mu$ equivalently defined as
\[
\M_\mu f(u,v,t)~=~
\sup_{\Rec\ni(u,v,t)}\vol\{\Rec\}^{-1}\iiint_{\Rec}\left|f(\xi,\eta,\tau+\mu(u\cdot\eta-v\cdot\xi))\right|d\xi d\eta d\tau.
\]
Observe that $\M_\mu$ is a special case of $\M_\rho^\alpha$ defined in (\ref{M_rho})-(\ref{rho_i}) at $\alpha=0$ with $\tau=x_1$ and $\rho_1=\mu(u\cdot\eta-v\cdot\xi)$. The $\L^p$-boundedness of $\M_\mu$ is proved by Christ \cite{Michael Christ 2}. Thereby, the elegant work is done by using  a mixture of techniques developed previously by Ricci and Stein \cite{Ricci-Stein 1}-\cite{Ricci-Stein 2} and Christ \cite{Michael Christ 1} for singular integrals defined on sub-manifolds within a general setting of Nilpotent Lie groups. 

$\diamond$ {\small Throughout, $\C>0$ is regarded as a generic constant depending on its sub-indices.}

In this paper, we prove the $\L^p\mt\L^q$-boundedness of $\M^\alpha_\rho$ by applying a geometric covering lemma due to C\'{o}rdoba and Fefferman \cite{Cordoba-Fefferman}.

{\bf Covering lemma} ~~
{\it Let $\Rec_j, j=1,2,\ldots$ be a sequence of rectangles in $\R^\N$ parallel to the coordinates. Denote $\chi$ to be an indicator function. There exists a subsequence $\Hat{\Rec}_k, k=1,2,\ldots$ such that
\bel{EST1}
\vol\Bigg\{\Cup_j \Rec_j\Bigg\}~\lesssim~\vol\Bigg\{ \Cup_k \Hat{\Rec}_k\Bigg\}
\eeq
and
\bel{EST2}
\left\|\sum_k \chi_{\Hat{\Rec}_k}\right\|_{\L^p(\R^\N)}~\leq~\C_p~\vol\Bigg\{\Cup_k\Hat{\Rec}_k\Bigg\}^{1\over p},\qquad 1<p<\infty.
\eeq}\\
\begin{remark}
The covering lemma above is true for any absolutely continuous measure satisfying the $\A_\infty$-property uniformly on every coordinate subspace. It's proof relies on a beautiful 'slicing' idea which reduces our assertion to the ($\N-1$)-dimensional subspace. For the Lebesgue measure case, we can prove this result with a more direct approach. 
\end{remark}

{\bf Theorem One}~~{\it Let $\M_\rho^\alpha$ defined in (\ref{M_rho})-(\ref{rho_i}) for $0\leq\alpha<1$. We have
\bel{Result One}
\left\|\M_\rho^\alpha f\right\|_{\L^q(\R^\N)}~\leq~\C_p\left\| f\right\|_{\L^p(\R^\N)},\qquad \alpha~=~{1\over p}-{1\over q},\qquad1<p\leq q<\infty.
\eeq}
We prove {\bf Theorem One} in the next section. We give a simpler proof of {\bf Covering lemma} in section 3.

\section{Proof of Theorem One}
\setcounter{equation}{0}
Let $\M_\rho^\alpha$ defined in (\ref{M_rho})-(\ref{rho_i}). Note that 
$\rho_i=\rho_i(y_{i+1},\ldots,y_\N,x_{i+1},\ldots,x_\N)$, $i=1,2,\ldots,\N-1$.
We have
\bel{f L^p-norm}
\begin{array}{lr}\ds
\int_{\R^\N}\left|f(y_1+\rho_1, y_2+\rho_2,\ldots,y_{\N-1}+\rho_{\N-1}, y_\N)\right|^p dy
\\\\ \ds
=~ \idotsint_{\R^{\N-1}} \left\{\int_\R \left|f(y_1+\rho_1, y_2+\rho_2,\ldots,y_{\N-1}+\rho_{\N-1}, y_\N)\right|^pdy_1\right\}  dy_2\cdots dy_\N
\\\\ \ds
=~\idotsint_{\R^{\N-1}} \left\|f(\cdot, y_2+\rho_2,\ldots,y_{\N-1}+\rho_{\N-1}, y_\N)\right\|^p_{\L^p(\R)} dy_2\cdots dy_\N
\\ \ds~~~~~~~~~~~~~~~~~~~~~~~~~~~~~~~~~~~~~~~~~~~~~~~~~~~~
\vdots
\\ \ds
=~\int_\R \left\|f(\cdot,  y_\N)\right\|^p_{\L^p(\R^{\N-1})} dy_\N~=~\left\| f\right\|^p_{\L^p(\R^\N)},\qquad 1<p<\infty.
\end{array}
\eeq
Given $\lambda>0$, we consider
\bel{U_lambda}
\U_\lambda~=~\Bigg\{ x\in\R^\N\colon \M^\alpha_\rho f(x)>\lambda\Bigg\}.
\eeq
For every $x\in\U_\lambda$, there is a rectangle $\Rec_j\subset\R^\N$ containing $x$ such that
\bel{R_j est}
 \vol\{\Rec_j\}^{\alpha-1}\int_{\Rec_j}\left|f(y_1+\rho_1, y_2+\rho_2,\ldots,y_{\N-1}+\rho_{\N-1}, y_\N)\right| dy~>~{1\over 2} \lambda.
\eeq

Let $x$ run through $\U_\lambda$. We find
\bel{U cover}
\U_\lambda~\subset~\Cup_j~\Rec_j.
\eeq
Recall (\ref{EST1})-(\ref{EST2}) in the {\bf Covering lemma}. Let $\Hat{\Rec}_k, k=1,2,\ldots$ be the selected subsequence of $\Rec_j, j=1,2,\ldots$. We have
\bel{union size Est}
\begin{array}{lr}\ds
\vol\Bigg\{ \Cup_k \Hat{\Rec}_k\Bigg\}
~\leq~\sum_k\vol \left\{ \Hat{\Rec}_k\right\} 
\\\\ \ds~~~~~~~~~~~~~~~~~~~~~
~\leq~ \sum_k \left\{2\lambda^{-1}\int_{\Hat{\Rec}_k}\left|f(y_1+\rho_1, y_2+\rho_2,\ldots,y_{\N-1}+\rho_{\N-1}, y_\N)\right| dy\right\}^{1\over 1-\alpha}
\\\\ \ds~~~~~~~~~~~~~~~~~~~~~
~\lesssim~  \lambda^{-{1\over 1-\alpha}}\left\{\sum_k \int_{\Hat{\Rec}_k}\left|f(y_1+\rho_1, y_2+\rho_2,\ldots,y_{\N-1}+\rho_{\N-1}, y_\N)\right| dy\right\}^{1\over 1-\alpha}\qquad\hbox{\small{($0\leq\alpha<1$)}}
\\\\ \ds~~~~~~~~~~~~~~~~~~~~~
~=~\lambda^{-{1\over 1-\alpha}}\left\{\int_{\R^\N}\left|f(y_1+\rho_1, y_2+\rho_2,\ldots,y_{\N-1}+\rho_{\N-1}, y_\N)\right| \sum_k \chi_{\Hat{\Rec}_k}(y)dy\right\}^{1\over 1-\alpha}
\\\\ \ds~~~~~~~~~~~~~~~~~~~~~
\leq~\lambda^{-{1\over 1-\alpha}} \left\{ \int_{\R^\N}\left|f(y_1+\rho_1, y_2+\rho_2,\ldots,y_{\N-1}+\rho_{\N-1}, y_\N)\right|^p dy\right\}^{{1\over p}{1\over 1-\alpha}}\left\|\sum_k\chi_{\Hat{\Rec}_k}\right\|_{\L^{p\over p-1}(\R^\N)}^{1\over 1-\alpha}
\\ \ds~~~~~~~~~~~~~~~~~~~~~~~~~~~~~~~~~~~~~~~~~~~~~~~~~~~~~~~~~~~~~~~~~~~~~~~~~~~~~~~~~~~~~~~~~~~~~~~~~~~~~~~~~~~~~~~~~~
 \hbox{\small{ by H\"older inequality}}
 \\ \ds~~~~~~~~~~~~~~~~~~~~
~\leq~\C_p~\lambda^{-1{1\over 1-\alpha}}\left\| f\right\|_{\L^p(\R^\N)}\vol\Bigg\{ \Cup_k \Hat{\Rec}_k\Bigg\}^{{p-1\over p}{1\over 1-\alpha}}
\qquad \hbox{\small{ by (\ref{f L^p-norm}) and (\ref{EST2})}}.
\end{array}
\eeq
By raising both sides of (\ref{union size Est}) to the $(1-\alpha)$-th power and taking into account for  $1-\alpha-{p-1\over p}={1\over p}-\left[{1\over p}-{1\over q}\right]={1\over q}$, we find
\bel{union size Est'}
\vol\Bigg\{ \Cup_k \Hat{\Rec}_k\Bigg\}^{1\over q}~\leq~\C_p~ \lambda^{-1}\left\| f\right\|_{\L^p(\R^\N)}.
\eeq
Let $\U_\lambda$ defined in (\ref{U_lambda}).
We obtain
\bel{weak type (p,p)}
\begin{array}{lr}\ds
\vol\Bigg\{ x\in\R^\N\colon \M_\rho^\alpha f(x)>\lambda\Bigg\}^{1\over q} ~=~\vol\left\{ \U_\lambda\right\}^{1\over q}
\\\\ \ds~~~~~~~~~~~~~~~~~~~~~~~~~~~~~~~~~~~~~~~~~~~~~~
~\leq~\vol\Bigg\{ \Cup_j \Rec_j\Bigg\}^{1\over q}\qquad \hbox{\small{by (\ref{U cover})}}
\\\\ \ds~~~~~~~~~~~~~~~~~~~~~~~~~~~~~~~~~~~~~~~~~~~~~~
~\lesssim~\vol\Bigg\{ \Cup_k \Hat{\Rec}_k\Bigg\}^{1\over q} \qquad\hbox{\small{by (\ref{EST1})}}
\\\\ \ds~~~~~~~~~~~~~~~~~~~~~~~~~~~~~~~~~~~~~~~~~~~~~~
~\leq~\C_p~\lambda^{-1} \left\|f\right\|_{\L^p(\R^\N)}\qquad \hbox{\small{by (\ref{union size Est'}) }}.
\end{array}
\eeq
We finish the proof of {\bf Theorem One} by using the weak type ($p,q$)-estimate in (\ref{weak type (p,p)}) and applying Marcinkiewicz interpolation theorem.

\section{A simpler proof of the covering  lemma}
\setcounter{equation}{0}
Given $\Rec_j, j=1,2,\ldots$, we select $\Hat{\Rec}_k, k=1,2,\ldots$  as follows.

Let $\Hat{\Rec}_1=\Rec_1$. Having chosen $\Hat{\Rec}_1, \Hat{\Rec}_2,\ldots,\Hat{\Rec}_{k-1}$, we pick $\Hat{\Rec}_k$ as the first rectangle $\Rec$ on the list of $\Rec_j$'s after $\Hat{\Rec}_{N-1}$ so that
\bel{R condition}
\vol\left\{~\Rec\cap\bigcup_{\ell=1}^{k-1}\Hat{\Rec}_\ell~\right\}~\leq~\frac{1}{2}\vol\left\{\Rec\right\}.
\eeq
Suppose that $\Rec$ is an unselected rectangle. There is a positive number $N$ such that $\Rec$ is on the list of $\Rec_j$'s after $\Hat{\Rec}_N$ and
\bel{unselect}
\vol\left\{~\Rec\cap\bigcup_{k=1}^N \Hat{\Rec}_k~\right\}~>~\frac{1}{2}\vol\left\{\Rec\right\}.
\eeq
Let $\M$ be the strong maximal operator defined as
\bel{M}
\M f(x)~=~\sup_{\Rec\ni x} \vol\{\Rec\}^{-1}\int_\Rec |f(y)|dy.
\eeq
From (\ref{unselect}),  we conclude
\bel{M>1/2}
\M \chi_{\Cup_k \Hat{\Rec}_k} (x)~>~{1\over 2},\qquad x\in\Cup_j  \Rec_j.
\eeq
By using the $\L^p$-boundedness of $\M$, we have
\bel{bound}
\begin{array}{lr}\ds
\vol\Bigg\{\Cup_j\Rec_j \Bigg\}~=~\int_{\Cup_j \Rec_j } dx
\\\\ \ds~~~~~~~~~~~~~~~~~~~~~
~\leq~2^2 \int_{\Cup_j \Rec_j } \Big[\M \chi_{\Cup_k \Hat{\Rec}_k}\Big]^2 (x) dx \qquad\hbox{\small{by (\ref{M>1/2})}}
\\\\ \ds~~~~~~~~~~~~~~~~~~~~~
~\lesssim~ \int_{\R^\N } \Big[\M \chi_{\Cup_k \Hat{\Rec}_k}\Big]^2 (x) dx
\\\\ \ds~~~~~~~~~~~~~~~~~~~~~
~\lesssim~ \int_{\R^\N} \Big[ \chi_{\Cup_k \Hat{\Rec}_k}\Big]^2 (x) dx
\\\\ \ds~~~~~~~~~~~~~~~~~~~~~
~=~\int_{\Cup_k \Hat{\Rec}_k}dx~=~\vol\Bigg\{\Cup_k\Hat{\Rec}_k \Bigg\}.
\end{array}
\eeq
On the other hand, (\ref{R condition}) suggests
\bel{R condition <}
\vol\left\{~\Hat{\Rec}_k\cap\bigcup_{\ell=1}^{k-1}\Hat{\Rec}_\ell~\right\}~\leq~\frac{1}{2}\vol\left\{\Hat{\Rec}_k\right\},\qquad k=1,2,\ldots
\eeq
which further implies 
\bel{vol E>S}
\vol\Bigg\{ \Hat{\Rec}_k\setminus\Cup_{\ell=1}^{k-1}\Hat{\Rec}_\ell\Bigg\}~>~\frac{1}{2}\vol\left\{\Hat{\Rec}_k\right\},\qquad k=1,2,\ldots.
\eeq
Let $\phi\in\L^{p\over p-1}(\R^\N)$ and $\left\|\phi\right\|_{\L^{p\over p-1}(\R^\N)}=1$. 
We have
\bel{integration}
\begin{array}{lr}\ds
\int_{\R^n}\phi(x)\sum_k\chi_{\Hat{\Rec}_k}(x) dx~=~\sum_k\int_{\Hat{\Rec}_k}\phi(x) dx
\\\\ \ds~~~~~~~~~~~~~~~~~~~~~~~~~~~~~~~~~~~~
~=~\sum_k \left\{\vol\{\Hat{\Rec}_k\}^{-1}\int_{\Hat{\Rec}_k}\phi(x) dx\right\} \vol\left\{\Hat{\Rec}_k\right\}
\\\\ \ds~~~~~~~~~~~~~~~~~~~~~~~~~~~~~~~~~~~~
~<~2\sum_k \left\{\vol\{\Hat{\Rec}_k\}^{-1}\int_{\Hat{\Rec}_k}|\phi(x)| dx\right\} \vol\Bigg\{ \Hat{\Rec}_k\setminus\Cup_{\ell=1}^{k-1}\Hat{\Rec}_\ell\Bigg\}\qquad\hbox{\small{by (\ref{vol E>S})}}
\\\\ \ds~~~~~~~~~~~~~~~~~~~~~~~~~~~~~~~~~~~~
~\lesssim~\sum_k \int_{\Hat{\Rec}_k\setminus\Cup_{\ell=1}^{k-1}\Hat{\Rec}_\ell} \left\{\vol\{\Hat{\Rec}_k\}^{-1}\int_{\Hat{\Rec}_k}|\phi(x)| dx\right\} dy
\\\\ \ds~~~~~~~~~~~~~~~~~~~~~~~~~~~~~~~~~~~~
~\lesssim~\sum_k\int_{\Hat{\Rec}_k\setminus\Cup_{\ell=1}^{k-1}\Hat{\Rec}_\ell}\M\phi(y) dy
~=~\int_{\Cup_k\Hat{\Rec}_k}\M\phi(y) dy.
\end{array}
\eeq
By applying H\"{o}lder inequality and using the $\L^p$-boundedness of $\M$, we find
\bel{boundedness of M}
\begin{array}{lr}\ds
\int_{\Cup_k\Hat{\Rec}_k}\M\phi(x) dx~\leq~\left\|\M\phi\right\|_{\L^{p\over p-1}(\R^\N)}\vol\Bigg\{\Cup_k\Hat{\Rec}_k\Bigg\}^{1\over p}
\\\\ \ds~~~~~~~~~~~~~~~~~~~~~~~~~~
~\leq~\C_p~\left\|\phi\right\|_{\L^{p\over p-1}(\R^\N)}~\vol\Bigg\{\Cup_k\Hat{\Rec}_k\Bigg\}^{1\over p}
~=~\C_p~\vol\Bigg\{\Cup_k\Hat{\Rec}_k\Bigg\}^{1\over p}.
\end{array}
\eeq
By substituting (\ref{boundedness of M}) to (\ref{integration}) and taking the supremum of $\phi$, we arrive at
\bel{p norm}
\left\|\sum_k\chi_{\Hat{\Rec}_k}\right\|_{\L^p(\R^\N)}~\leq~\C_p~\vol\Bigg\{\Cup_k\Hat{\Rec}_k\Bigg\}^{\frac{1}{p}}.
\eeq

{\small  Westlake University, Hangzhou, 310030, China}\\
{\small email: wangzipeng@westlake.edu.cn}

\end{document}